\documentclass[leqno]{article}[12pt]
\usepackage{enumitem}
\usepackage{mathrsfs}
\usepackage{geometry}
\usepackage[toc,page]{appendix}
\usepackage{hyperref}
\usepackage{stmaryrd}
\hypersetup{colorlinks=false,  linkcolor=black}
\usepackage{amsfonts}
\usepackage{amssymb}
\usepackage{ae,aecompl}
\usepackage{float}
\usepackage[nameinlink,capitalize]{cleveref}
\usepackage{amsmath}
\usepackage{amsthm}
\usepackage{graphicx}
\usepackage{enumitem}
\usepackage{epstopdf}
\usepackage{microtype}
\usepackage{algorithmic}
\usepackage{mathrsfs}
\usepackage{esint}
\usepackage{url}
\usepackage{mathtools}
\usepackage{arydshln}
\usepackage{chngcntr}
\usepackage[utf8]{inputenc}
\usepackage{xspace}
\usepackage{setspace}
\usepackage{nicefrac}
\usepackage{xcolor}
\usepackage[style=alphabetic, firstinits=true, sorting=debug]{biblatex}
\numberwithin{equation}{section}

\theoremstyle{plain}
\newtheorem{lem}{Lemma}[section]

\newtheorem*{thm*}{Theorem}

\newtheorem*{conj*}{Conjecture}
\theoremstyle{definition}

\newtheorem*{rmks*}{Remarks}

\newtheorem*{rmk*}{Remark}
\newtheorem*{swc}{S-W condition}

\newcommand{\select}{ \ ; \ }

\newcommand{\Z}{\mathbb{Z}}

\newcommand{\R}{\mathbb{R}}

\newcommand{\ov}{\overline}
\newcommand{\sumsum}{\mathop{\sum \sum}}
\newcommand{\sump}{\sideset{}{'}\sum}
\newcommand{\sumst}{\sideset{}{^*}\sum}
\newcommand{\dum}{\, d}

\DeclareMathOperator{\im}{Im}

\makeatletter
\let\@@pmod\pmod
\DeclareRobustCommand{\pmod}{\@ifstar\@pmods\@@pmod}
\def\@pmods#1{\mkern4mu({\operator@font mod}\mkern 6mu#1)}
\makeatother

\DeclareFieldFormat{labelalpha}{\thefield{entrykey}}

\begin{document}


\title{Sums over Vanishing Determinants}
\author{John Friedlander, Henryk Iwaniec}
\date{}
\maketitle
\tableofcontents


\section{Introduction}

Among the simplest determinants are those which appear in the resolution of systems of two linear equations: 
\begin{equation}
\label{1.1}
\begin{split}
  u_1 x + v_1 y &= a , \\
  u_2 x + v_2 y &= b .
\end{split}
\end{equation}
We are interested in the system \eqref{1.1} over the rational integers $\Z$  as well as over the ring $\Z / q Z$. If 
\[ \det \left( \begin{array}{cc} u_1 v_1 \\ u_2 v_2 \end{array} \right) = u_1 v_2 - u_2 v_1 \neq 0, \]
then the system has a unique solution, not necessarily integral. Otherwise, the linear space of solutions, if indeed there are any, has dimension one. Despite the relative simplicity of these $2 \times 2$ determinants, they arise in numerous arithmetic problems. We found our motivation from specific applications, for example its relevance to the problem studied in [FI2]. Therein, questions arise about the summation of arithmetic functions, such as the M\"{o}bius function, over singular systems modulo $q$. This means that we encounter sums of the type
\begin{equation}
  S_q = \sum_{u_1 v_2 \equiv u_2 v_1 \pmod* q} f(u_1 , v_1) g(u_2 , v_2).
  \label{1.2}
\end{equation}
One may say we are summing over the vectors $[u_1, v_1]$ and $[v_2, - u_2]$ which are orthogonal modulo $q$.

For notational convenience we prefer to use arithmetic functions $c(z)$ defined on Gaussian integers $z = u + iv$; in fact the complex analytic structure of $c(z)$ will be present. Let us assume for simplicity that $f$, $g$ are complex conjugate so our determinant sums are
\begin{equation}
  S_q = \sumsum_{u_1 v_2 \equiv u_2 v_1 \pmod* q} c(z_1) \overline{c}(z_2)
  = \sumsum_{\im z_1 \overline{z}_2 \equiv 0 \pmod* q} c(z_1) \overline{c}(z_2).
  \label{1.3}
\end{equation}
Furthermore, for a technical reason, we assume that $c(z)$ is supported on odd primitive Gaussian integers $z \neq \pm 1, \pm i$, which means $z$ has coordinates of different parity, coprime and non-zero:
\begin{equation}
  z = u + iv, \quad u \not\equiv v ({\rm mod}\, 2), \quad (u,v) = 1,\  uv \neq 0.
  \label{1.4}
\end{equation}
Note that $(z, \overline{z}) = 1$.

\smallskip

Our main result gives a non-trivial bound for the sum of $ S_q$ over $q\leqslant Q$ provided 
  $c(z)$ is supported on odd primitive numbers $1 < |z|^2 \leqslant X$, $|c(z)| \leqslant 1$, and that it satisfies a condition (of Siegel-Walfisz type) for uniformity of distribution in arithmetic progressions to small moduli. In this case we are able to prove a non-trivial bound with $q$ running, within a few logarithms, up to $X$. See Section 3 for the precise statement and further remarks. 

\smallskip

If one were to restrict $z$ to Gaussian integers coprime with $q$, then the congruence $\im z_1 \overline{z}_2 \equiv 0 \pmod q$ could be decoupled; it is equivalent to
\begin{equation}
  \frac{\overline{z}_1}{z_1} \equiv \frac{\overline{z}_2}{z_2} ({\rm mod}\, q).
  \label{1.5}
\end{equation}
This congruence can be detected by the orthogonality of characters $\psi ({\rm mod}\, q)$ on the group
\begin{equation}
  G = \left\{ \alpha \in \Z[i] / q \Z[i] \select \alpha \overline{\alpha} \equiv 1 ({\rm mod}\, q) \right\}.
  \label{1.6}
\end{equation}
We obtain
\begin{equation}
\label{1.7}
\sumsum_{\substack{\im z_1 \ov{z}_2 \equiv 0\, ({\rm mod}\, q) \\ (z_1 z_2, q) = 1 }} c(z_1) \ov{c}(z_2) 
  = \frac{1}{|G|} \sum_{\psi \in \hat{G}} \left| \sum_{(z,q) = 1} c(z) \psi(\ov{z})/ \psi(z) \right|^2. 
\end{equation}
Note that for $q$ squarefree
\begin{equation}
  |G| = \prod_{p | q} (p - \chi_4(p)),
  \label{1.8}
\end{equation}
while the full group of classes $\alpha ({\rm mod}\, q)$, $(\alpha, q) = 1$ has larger order, namely $\varphi(q)|G|$.

\smallskip

Although this expression of the determinant 
sums in terms of characters presents quite an attractive option,
we are going to work with congruences modulo $q$ for greater transparency. Actually, this turns out to be necessary for the performance of certain transformations, such as switching moduli and enlarging moduli. These techniques, carried out here in Sections 5 and 6 respectively, are reminiscent of the strategy used previously in Sections 12 and 13 of our earlier work [FI1].

\smallskip

The partially reduced determinant sum
\begin{equation}
  S'_q = \sumsum_{\substack{u_1 v_2 \equiv u_2 v_1 \pmod* q \\ (v_1 v_2, q) = 1}} c(z_1) \ov{c}(z_2)
  \label{1.9}
\end{equation}
is not very different from \eqref{1.3}. Precisely, the original sum decomposes as
\begin{equation}
  S_q = \sum_{d |q} S'_{q/d}
  \label{1.10}
\end{equation}
where the reduced sum $S'_{q/d}$ has coefficients $c_d(u+iw) = c(u+idw)$. The large divisors $d$ can be easily eliminated up to negligible error terms and the small ones make an insignificant deformation of the coefficients. Now, the congruence in the reduced sum can be decoupled as
\begin{equation}
  u_1 \ov{v}_1 \equiv u_2 \ov{v}_2 \pmod q
  \label{1.11}
\end{equation}
where $\ov{v}$ represents the multiplicative inverse of $v$ modulo $q$ and not complex number conjugation. Detecting \eqref{1.11} by additive characters we get
\begin{equation}
 0  \leqslant  S'_q = \frac{1}{q} \sum_{a \pmod* q} \left| \sump_{z = u + iv} c(z) e\left( \frac{a}{q} u \ov{v} \right) \right|^2   \leqslant S_q ,
  \label{1.12}
\end{equation}
where the latter inequality comes from positivity and the observation of the term $d=1$ in \eqref{1.10}. 
Here and thereafter, the stroke restricts the summation to Gaussian integers 
$z$ with $(q, \im z) = 1$. 
 
\smallskip

If the coefficients can be decoupled along the coordinates, say
\begin{equation}
  c(z) = c(u+iv) = a(u) b(v),
  \label{1.13}
\end{equation}
then we have a more favorable expression
\begin{equation}
  S'_q = \frac{1}{q} \sum_{a \pmod* q} \left| \sum_{\substack{z = u+iv \\ (v,q) =1}} a(u) \ov{b}(v) e\left( \frac{a}{q} u v \right) \right|^2.
  \label{1.14}
\end{equation}
To see this, open the square and change $a$ into $a v_1 v_2$ modulo $q$. Here of course, $\ov{b}(v)$ does stand for the complex conjugate of $b(v)$. This expression can be treated directly by the classical large sieve inequality (consider $n = uv$ as a single variable) producing very strong estimates. Unfortunately, \eqref{1.13} does not hold in the most interesting cases and we shall need to proceed without use of this property, a task which constitutes the core of our work.

\section{A few observations about other restricted sums}

%

When restricted to numbers $z = u + iv$ with $(uv, q) =1$,
our determinant sum modulo $q$, now denoted by $S''_q$, enjoys the symmetry
\begin{equation}
\label{2.1}
  S''_q = \frac{1}{q} \sum_{a \pmod* q} \left| \sum_{\substack{z = u + iv \\ (uv,q) = 1}} c(z) e\left( \frac{a}{q} u \ov{v} \right) \right|^2 
 =  \frac{1}{q} \sum_{a \pmod* q} \left| \sum_{\substack{z = u + iv \\ (uv,q) = 1}} c(z) e\left( \frac{a}{q} \ov{u} v \right) \right|^2\ . 
\end{equation}
To see this open the square and change $a$ into $a \ov{u}_1 \ov{u}_2 v_1 v_2$ modulo $q$.

\smallskip

The determinant sum modulo $q$ which is restricted to numbers $z = u + iv$ with $(v,q) = 1$ and $(z,q) = 1$, denoted by $S_q'''$, satisfies 
\begin{equation}
 S_q''' = \frac{1}{q} \sum_{a \pmod* q}  \left| \sum_{\substack{z = u + iv \\ (v,q) = (z,q) = 1}} c(z) e\left( \frac{a}{q} u \ov{v} \right) \right|^2 \leqslant S_q'\ .
  \label{2.3}
\end{equation}
To see it open the square and write the congruence $\im z_1 \ov{z}_2 \equiv 0 \pmod q$ in the form $z_1 \ov{z}_2 \equiv \ov{z}_1 z_2 ({\rm mod}\, 2q)$. Hence, one of the two coprimality conditions $(z_1, q) = 1$, $(z_2,q) = 1$, is redundant and
\begin{equation}
S'''_q = \frac{1}{q} \sum_{a \pmod* q} \left( \sump_{(z,q) = 1} c(z) e\left( \frac{a}{q} u \ov{v} \right) \right) \left( \sump_{z} \ov{c}(z) e\left( -\frac{a}{q} u \ov{v} \right) \right).
 \label{2.4}
\end{equation}
Applying Cauchy's inequality we derive $S'''_q \leqslant S'_q$.

\section{Statement of the theorem}

Let $S_q(X)$ and $S'_q(X)$ denote our determinant sums modulo $q$ with the coefficients $c(z)$ cropped to the disc
\begin{equation}
  |z|^2 \leqslant X . 
  \label{3.1}
\end{equation}
Our goal is to estimate these sums on average over $q \leqslant Q$, that is the sums
\begin{align}
  S(Q,X) &= \sum_{q \leqslant Q} S_q (X) \label{3.2}= \sum_{q\leqslant Q}\, \sumsum_{u_1 v_2 \equiv u_2 v_1 \,({\rm mod}\,q) } c(u_1+iv_1) \overline{c}(u_2+iv_2) , \\
  S'(Q,X) &= \sum_{q \leqslant Q} S'_q(X)= \sum_{q \leqslant Q} \frac{1}{q} \sum_{a \pmod* q} \left| \sump_{\substack{z = u + iv\\ (v,q)=1}} c(z) e\left( \frac{a}{q} u \ov{v} \right) \right|^2 \label{3.3}\ .
\end{align}
Note that
\begin{equation}
 0 \leqslant S'(Q,X)\leqslant S(Q,X)\ .
  \label{3.4}
\end{equation}

We assume that
\begin{equation}
  |c(z)| \leqslant 1.
  \label{3.5}
\end{equation}
Then the trivial estimation of \eqref{1.3} yields
\begin{equation}
  S(Q,X) \ll (Q + X L^5)X, \quad L = \log X,
  \label{3.6}
\end{equation}
and we wish to improve it by a factor of an arbitrary power of $L$ with $Q$ nearly as large as $X$. The saving factor need not be very large, yet it is crucial for applications. But, note that the residue class $a \pmod q$ in \eqref{1.12} need not be reduced so in general \eqref{3.6} cannot be improved. However, it is possible to beat \eqref{3.6} if the sequence of coefficients $c(z)$ admits a considerable cancellation in sums over residue classes in the Gaussian domain to 
small moduli. We assume the following:
\begin{swc}
  Let $k \geqslant 0$, $\ell \geqslant 1$, $\alpha \in \Z[i]$ and $t \in \R$. Then, we have
  \begin{equation}
    \sum_{z \equiv \alpha \pmod* \ell} c(z) (\ov{z} / z)^k |z|^{it} \ll \ell^2 (k^2 + 1)(t^2 + 1) X L^{-B} ,
    \label{3.7}
  \end{equation}
  where $B$ is any positive number and the implied constant depends only on $B$.
\end{swc}
Our main example of $c(z)$ satisfying \eqref{3.7} is the M\"{o}bius function
\begin{equation}
  c(z) = \mu(|z|^2).
  \label{3.8}
\end{equation}
In this case, \eqref{3.7} is just the Siegel-Walfisz condition in the Gaussian domain.
\begin{thm*}
  Suppose $c(z)$ is supported on odd primitive numbers $1 < |z|^2 \leqslant X$, $|c(z)| \leqslant 1$, and that it satisfies \eqref{3.7}. Then, we have
  \begin{equation}
    S(Q,X) = \sum_{q\leqslant Q}\, \sumsum_{u_1 v_2 \equiv u_2 v_1 \,({\rm mod}\,q) } c(u_1+iv_1) \overline{c}(u_2+iv_2)
\ll \left(QL^{28 A + 148} + XL^{-A}\right)X , 
    \label{3.9}
  \end{equation}
  where $A$ is any positive number 
and the implied constant depends only on $A$.
\end{thm*}

\section{First estimation of $S(Q,X)$}
We begin with (see \eqref{1.12})
\begin{equation}
  S'(Q,X) = \sum_{q \leqslant Q} \frac{1}{q} \sum_{a \pmod* q} \left| \sump_{\substack{z = u + iv\\ (v,q)=1}} c(z) e\left( \frac{a}{q} u \ov{v} \right) \right|^2.
  \label{4.1}
\end{equation}
We write the fraction $a/q$ in its lowest terms and decompose \eqref{4.1} into $T_1 + T_2$, where
\[ T_1 = \sumsum_{\substack{qr \leqslant Q \\ q > Q_0}} \frac{1}{qr} \sumst_{a\, ({\rm mod}\, q)} \left| \sum_{\substack{z = u + iv \\ (v,qr) = 1}} c(z) e\left( \frac{a}{q} u \ov{v} \right) \right|^2. \]
and $T_2$ is the complementary sum in which $q \leqslant Q_0$. We estimate $T_1$ by the large sieve inequality as follows:
\begin{align*}
  T_1 &\leqslant \sqrt{X} \sum_v \sumsum_{\substack{qr \leqslant Q \\ q > Q_0}} \frac{1}{qr} \sumst_{a \pmod* q}  \left| \sum_{z = u + iv} c(z) e\left( \frac{a}{q} u \right) \right|^2 \\
  &\ll \sqrt{X} \sum_{r \leqslant Q} \frac{1}{r} \left( \frac{Q}{r} + \frac{\sqrt{X}}{Q_0} \right) X \ll \left( Q \sqrt{X} + \frac{X}{Q_0} \log Q\right) X.
\end{align*}
For estimation of $T_2$ we appeal to the S-W condition \eqref{3.7}. To this end we generalize \eqref{3.7} as follows:
\begin{swc}
  Let $e \geqslant 1$, $k \geqslant 0$, $\ell \geqslant 1$, $\alpha \in \Z[i]$ and $t \in R$. Then we have
  \begin{equation}
    \sum_{\substack{z = u + iv \equiv \alpha ({\rm mod}\,\ell) \\ (v,e) = 1}} c(z) (\ov{z}/z)^k |z|^{it} \ll \tau(e) \ell^2 (k^2 + 1)(t^2+1) X L^{-B} , 
    \label{4.2}
  \end{equation}
  where $B$ is any positive number and the implied constant depends only on $B$.

  \begin{proof}
    Relax the coprimality condition $(v,e) = 1$ by the M\"{o}bius formula. Accordingly the sum \eqref{4.2} splits into
    \[ \sum_{d|e} \mu(d) \sum_{\substack{z=u+iv\equiv \alpha ({\rm mod}\, \ell)\\ d|v}} c(z) (\overline z /z)^k|z|^{it}\ .\]
  For $d \leqslant L^B$ apply \eqref{3.7} with $\ell$ replaced by $d\ell$ and $B$ replaced by $3B$, getting the bound \eqref{4.2}. For $d > L^B$ we estimate trivially by $\tau(e) XL^{-B}$. \endproof 

\smallskip

  Now we are ready to apply \eqref{4.2} for the estimation of $T_2$. To this end fix $u$, $v$ modulo $q$ and apply \eqref{4.2} with $e = r$, $k= 0$, $\ell=q$ and $t=0$. The number of relevant classes $z \equiv \alpha \pmod q$ is $q^2$, so the inner sum over $z$ in $T_2$ is bounded by $\tau(r) q^4 X L^{-6B}$. Hence 
  \[ T_2 \ll \sumsum_{\substack{qr \leqslant Q \\ q \leqslant Q_0}} \frac{1}{qr} q\left( \tau(r) q^4 X L^{-6B} \right)^2 \ll Q_0^9 X^2 L^{4 -12 B}. \]
  Choosing $Q_0 = L^{B+1}$ and adding the above bounds for $T_1$ and $T_2$ we conclude that 
  \begin{equation}
    S'(Q,X) \ll (Q\sqrt{X} + X L^{-B}) X
    \label{4.3}
  \end{equation}
  where $B$ is any positive number and the implied constant depends only on $B$.
  \end{proof}
\end{swc}

  Next we use \eqref{4.3} to estimate $S(Q,X)$. We estimate $S'_{q/d}$ in \eqref{1.10} trivially by
  \begin{align*}
    \mathop{\sum \sum \sum \sum}_{\substack{u_1 w_2 \equiv u_2 w_1 \ ({\rm mod}\, q/d)  \\ (u_1, w_1) = (u_2, w_2) = 1 \\ u_1,u_2 < \sqrt{X} ;\, w_1, w_2 < \sqrt{X} /d}} 1 & \ll \frac{X}{d} + \sum_{1 \leqslant n \leqslant X/q} \sumsum_{\substack{1 \leqslant a, b < X/d \\ |a \pm b| = nq / d}} \tau(a) \tau(b) \\
    &\ll \frac{X}{d} + \sum_{1 \leqslant n < X/q} \sum_{1 \leqslant a < X/d} \tau(a)^2 \ll \left( 1 + \frac{X}{q} L^3 \right) \frac{X}{d}.
  \end{align*}
  Hence
  \[ S(Q,X) = \sum_{q \leqslant Q} \sum_{\substack{d | q \\ d \leqslant L^C}} S'_{q/d} (X) + O\bigl((Q + XL^5) X L^{-C}\bigr)\ , \]
where $C\geqslant 1$ is a constant at our disposal. 
  This gives
  \begin{equation}
    S(Q,X) = \sum_{d \leqslant L^C} S'(Q/d,X) + O\bigl( (Q + XL^5) XL^{-C}\bigr) ,
    \label{4.4}
  \end{equation}
  where the coefficients in $S'(Q/d, X)$ are $c_d(u+iv) = c(u+idv)$. These coefficients satisfy the condition \eqref{4.2} with $\ell$ replaced by $d\ell$; hence the upper bound \eqref{4.2} is larger by a factor $d^2 \leqslant L^{2C}$. 
This factor can be ignored because $B$ in \eqref{4.2} is arbitrary. Therefore \eqref{4.3} applies to every $S'(Q/d, X)$ in \eqref{4.4} giving
  \[ S(Q,X) \ll (Q \sqrt{X} L^{2C + 1} + X L^{-B} )X + (Q + X L^5)XL^{-C}. \]
Hence we conclude (choose $B=A$ and $C = A + 5$):
\begin{lem}
  For $Q \geqslant 1$, $X \geqslant 2$ we have
  \begin{equation}
    S(Q,X) \ll \left(Q \sqrt{X} L^{2A + 10} + X L^{-A}\right) X
    \label{4.5}
  \end{equation}
  where $A$ is any positive number and the implied constant depends only on $A$.
\end{lem}
Note that Lemma 4.1 yields \eqref{3.9} if $Q\leqslant \sqrt X L^{-3A-10}$ . 

\section{Second estimation of $S(Q,X)$}

In this section we reduce the problem for larger moduli to that for smaller ones by switching divisors. For a technical reason we subdivide the range of moduli into dyadic segments. Since $S_q (X)\geqslant 0$ it suffices to estimate the weighted sum
\begin{equation}
  S(\sim Q, X) = \sum_q W\left( \frac{q}{Q} \right) S_q (X)
  \label{5.1}
\end{equation}
where $W(x) \geqslant 0$ is a fixed smooth function supported on $1 < x < 4$ with $W(2) \geqslant 1$. Recall that
\begin{equation}
  S_q (X) = \sum_{z_1 \ov{z}_2  \equiv \ov{z}_1 z_2 \pmod* q} c(z_1) \ov{c} (z_2)
  \label{5.2}
\end{equation}
and $c(z)$ are cropped to the disc $|z|^2 \leqslant X$.

The determinant $z_1 \ov{z}_2 - \ov{z}_1 z_2$ vanishes only for the diagonal 
terms $z_1 = z_2$ which yields the contribution
\begin{equation}
  V_0(Q,X) \ll QX.
  \label{5.3}
\end{equation}
On the off-diagonal we have
\begin{equation}
  |z_1 \ov{z}_2 - \ov{z}_1 z_2| = dq \,\, \text{ with }\,\, 1 \leqslant d \leqslant 2XQ^{-1}.
  \label{5.4}
\end{equation}
Let $C\geqslant 1$ be a constant to be specified later. If $d \leqslant X Q^{-1} L^{-C}$ we estimate trivially getting the contribution
\begin{equation}
  V_1(Q,X) \ll X^2 L^{4-C}.
  \label{5.5}
\end{equation}
We are left with
\begin{equation}
  V_2(Q,X) = \sum_d \sumsum_{z_1 \ov{z}_2 \equiv \ov{z}_1 z_2 ({\rm mod}\, d)} c(z_1) \ov{c}(z_2) W(|z_1 \ov{z}_2 - \ov{z}_1 z_2 |/dQ)
  \label{5.6}
\end{equation}
where $d$ runs over the segment
\begin{equation}
  DL^{-C} < d < 2D, \ \ D = X/Q.
  \label{5.7}
\end{equation}

We need to separate the variables $z_1$, $z_2$ involved in the weight function without contaminating the coefficients $c(z_1)$, $c(z_2)$ too much. For this job it is convenient to choose $W(x)$ in the form of the convolution
\begin{equation}
  W(x) = \int A(x/y) B(y) y^{-1} \dum y
  \label{5.8}
\end{equation}
where $A(x)$, $B(y)$ are supported on the segments $1 \leqslant x, y \leqslant 2$. We have
\begin{equation}
  |z_1 \ov{z}_2 - \ov{z}_1 z_2 | = |z_1 z_2| |z - 1|, \ \ z = z_1 \ov{z}_2 / \ov{z}_1 z_2,
  \label{5.9}
\end{equation}
and 
\begin{equation}
  W(|z_1 \ov{z}_2 - \ov{z}_1 z_2 |/dQ) = \int A \left( \frac{2|z_1 z_2|}{dQw} \right) B\left( \frac{w}{2} |z - 1| \right) \frac{dw}{w}.
  \label{5.10}
\end{equation}

It follows from the support of $A(x)$ that $dQw < 2|z_1 z_2| \leqslant 2X$, hence $w < 2 L^C$ by \eqref{5.7}, and it follows from the support of $B(y)$ that $1 \leqslant \frac{w}{2} | z - 1| \leqslant w$. Therefore the variable $w$ in the integral representation \eqref{5.10} runs over the segment
\begin{equation}
  1 \leqslant w \leqslant W = 2L^C\ .
  \label{5.11}
\end{equation}
Truncating the Mellin integral 
\[ A(x) = \int_{-\infty}^{\infty} \tilde{A}(t)x^{it} \dum t \]
to the segment $|t|\leqslant T=L^C$  
and using the bound $  \tilde{A}(t) \ll (t^2 + 1)^{-1} $, 
we obtain 
\begin{equation}
  A(x) = \int_{-T}^T \tilde{A}(t) x^{it} \dum t + O\left( T^{-1} \right).
  \label{5.12}
\end{equation}
The contribution to $V_2(Q,X)$ of the above error term is 
\begin{equation}
  V_3(Q,X) \ll X^2 L^{4 - C}
  \label{5.13}
\end{equation}
by a trivial estimation. Hence
\begin{equation}
  A \left( \frac{2|z_1 z_2|}{dQw} \right) = \int_{-T}^T \tilde{A}(t) \left( \frac{2 |z_1 z_2|}{dQw} \right)^{it} \dum t + O\left( T^{-1} \right).
\end{equation}

We treat $B\left( \frac{w}{2} |z -1| \right)$ by its Fourier series expansion. Putting $z = e^{2i\alpha}$ we write
\[ B\left( \frac{w}{2} |z-1| \right) = B(w |\sin \alpha|).\]
This is an even periodic function of $\alpha$ of period $\pi$ so we have
\[ B(w |\sin \alpha|) = \sum_0^\infty b_k(w) \cos(2\alpha k) \]
with the coefficients
\[  b_k(w) = \frac{1}{\pi}\int_0^\pi B(w \sin \alpha) \cos (2 \alpha k) \dum \alpha . \]
By the support of $B(y)$ it follows that $1 < w \sin \alpha < 2$. Hence the trivial estimation yields $b_k (w) \ll w^{-1}$. If $k > 0$ we can integrate by parts two times getting
\begin{align*}
  4 \pi k^2 b_k (w) &= - \int_0^\pi (B(w \sin \alpha))'' \cos (2\alpha k) \dum \alpha \\
  &= -w \int_0^\pi (B'(w \sin \alpha) \cos \alpha)' \cos (2 \alpha k) \dum \alpha.
\end{align*}
Hence the trivial estimation yields $b_k(w) \ll w k^{-2}$. Combining the two estimations we obtain
\[ b_k(w) \ll w(k + w)^{-2}. \]
Hence
\begin{equation}
  B\left( \frac{w}{2} |z-1| \right) = \sum_{0\leqslant k < K} b_k(w) \left( z^k + z^{-k} \right) + O(w K^{-1})
  \label{5.15}
\end{equation}
with any $K\geqslant 1$. 
We choose $K = L^{2C}$ so the error term contributes to $V_2(Q,X)$
\begin{equation}
  V_4(Q,X) \ll X^2 L^{4 -C}
  \label{5.16}
\end{equation}
by a trivial estimation. We are left with
\begin{equation}
  V_5(Q,X) = \sum_{|k| < K} \int_{-T}^T \tilde{A}(t)\left( \frac{2}{Q} 
\right)^{it}\nu(k,t)  S^{(kt)} (D,X) \dum t
  \label{5.17}
\end{equation}
where
\begin{equation}
  \nu(k,t) = \int_1^W b_k(w) w^{-1 - it} dw \ll 1,
  \label{5.18}
\end{equation}
and
\begin{equation}
  S^{(kt)} (D,X) = \sum_d \sum_{z_1\ov{z}_2\equiv\ov{z}_1z_2 ({\rm mod}\, d)} c(z_1) \ov{c}(z_2) |z_1 z_2|^{it} \left( \frac{\ov{z}_1}{z_1} \right)^k \left( \frac{\ov{z}_2}{z_2} \right)^{-k}.
  \label{5.19}
\end{equation}
Recall that $K = L^{2C}$, $T = L^C$, $W = 2L^C$ and $d$ runs over the segment \eqref{5.7}.

\smallskip

Our first estimation \eqref{4.5} is applicable to $S^{(kt)} (D,X)$ with $Q$ replaced by $2D = 2XQ^{-1}$ giving
\[ S^{(kt)}(D,X) \ll \left( D \sqrt{X} L^{2B + 10} + XL^{-B} \right) X. \]
Hence
\begin{equation}
  V_5 (Q,X) \ll L^{3C} \left(Q^{-1} X^{\frac{3}{2}} L^{2B + 10} + XL^{-B}\right)X.
  \label{5.20}
\end{equation}
Adding the bounds \eqref{5.3}, \eqref{5.5}, \eqref{5.13}, \eqref{5.16} to \eqref{5.20} (with $B$ changed into $B+3C$, as we may since $B$ is arbitrary),  we obtain
\[ S(\sim Q, X) \ll QX + X^2 L^{4 - C} + L^{3C} \left( Q^{-1}X^{\frac{3}{2}} L^{2B + 6C + 10} + XL^{-B - 3C} \right) X . \]
 Now, choosing $C = B+4$, we get
\begin{lem}
The weighted sum \eqref{5.1} satisfies
\begin{equation}
  S(\sim Q, X) \ll \left( Q + Q^{-1}X^{\frac{3}{2}} L^{11B + 46} + XL^{-B} \right) X ,
  \label{5.21}
\end{equation}
where $B$ is any positive number and the implied constant depends only on $B$.
\end{lem}

Note that, on taking $B=A+1$, \eqref{5.21} yields 
\begin{equation}
  S(\sim Q, X) \ll \left( Q + XL^{-A-1}\right)X
  \label{5.22}
\end{equation}
if $Q \geqslant \sqrt{X} L^{12A+ 47}$. We still need to cover the middle range 
\begin{equation}
\sqrt{X} L^{-3A -10} < Q < \sqrt{X} L^{12A + 47}.
  \label{5.23}
\end{equation}

\section{Leapfrog and completing the proof}

We are able to estimate the sum over this middle range of moduli by transforming it to a sum over larger moduli covered in the previous section. We accomplish this by enlarging $q$ artificially with the aid of the primes in a segment $P < p \leqslant 2P$, where $P$ is at our disposal. To simplify the notation we hide $X$ which controls the support of $c(z)$ and so we put $ S'_q =  S'_q (X)$. 
Given a prime $p$, we write
\[ S'_q = \sumsum_{\substack{u_1 v_2 \equiv u_2 v_1 \pmod* q \\ (v_1 v_2, pq) = 1 }} c(u_1 + i v_1) \ov{c}(u_2 + iv_2) + \sumsum_{\substack{u_1 w_2 \equiv u_2 w_1 \pmod* q \\ (p w_1 w_2, q) = 1}} c(u_1 + i p w_1) \ov{c}(u_2 + i p w_2). \]
Hence
\begin{equation*}
  S'_q  \leqslant  \frac{2}{q} \sum_{a \pmod* q} \left| \sum_{\substack{z = u+iv \\ (v,pq) =1}} c(z) e\left( \frac{a}{q} u \ov{v} \right) \right|^2 
  + \frac{2}{q} \sum_{a \pmod* q} \left| \sum_{\substack{z = u + ipw \\ (w,q) = 1}} c(z) e\left( \frac{a}{q} u \ov{w} \right) \right|^2.
\end{equation*}
Note that in the second sum we dropped the condition $(p,q) = 1$, as we can by positivity. The second sum is just the determinant sum modulo $q$ with the coefficients $c_p(u + iw) = c(u + ipw)$, which we denote by $S^{(p)}_q$. For large $p$ it will be sufficient to estimate $S^{(p)}_q$ trivially. In the first sum we write $a/q = pa/pq$ and ignore that $pa$ runs over the multiples of $p$, as we can by positivity. Hence the first sum is bounded by $p S'_{pq}$. We end up with 
the ``enlarging moduli inequality''
\begin{equation}
  S'_q \leqslant 2p S'_{pq} + 2 S^{(p)}_q\ ,
  \label{6.1}
\end{equation}
which holds for every $q$ and prime $p$.

\smallskip

The sum $S_q^{(p)}$ in \eqref{6.1} can be estimated trivially as follows: 
\begin{align*}
  S_q^{(p)} &\leqslant 4 \sumsum_{\substack{u_1 w_2 \equiv u_2 w_1 \pmod* q \\ (u_1, w_1) = (u_2, w_2) = 1 \\ 1 \leqslant u_1, u_2 < \sqrt{X} ;\, 1 \leqslant w_1, w_2 < \sqrt{X} /p}} 1 \\
  &\leqslant \frac{4X}{p} + 8 \sum_{1 \leqslant n \leqslant X/pq} \sum_{\substack{1 \leqslant a, b < X/p \\ |a \pm b| = nq}} \tau(a) \tau(b) \ll \left( 1 + \frac{X}{pq} L^3 \right) \frac{X}{p}.
\end{align*}
Put
\[ S'(\sim Q) = \sum_{Q< q \leqslant 2Q} S'_q. \]
Multiply \eqref{6.1} by $\log p$ and sum over $P< p \leqslant 2P$, $Q <q \leqslant 2Q$ getting:
\begin{lem}
  For any $P \geqslant 4$ and $Q \geqslant 1$ we have
  \begin{equation}
    S'(\sim Q) \ll \left(S'(\sim PQ) + S'(\sim 2PQ)\right) \log PQ + \left( Q + \frac{X}{P} L^4 \right) \frac{X}{P} , 
    \label{6.2}
  \end{equation}
  where the implied constant is absolute.
\end{lem}

From $S'(\sim Q)$ we go to $S(\sim Q,X)$ using the approximate formula \eqref{4.4}. First, for every $S'(\sim Q/d)$ in \eqref{4.4} we apply \eqref{6.2} getting
\begin{equation*}
  \begin{split}
  S(\sim Q,X)\ll \sum_{d \leqslant L^C} \left( S'(\sim PQ/d) + S'(\sim 2PQ/d) \right) \log PQ \\
  + \left( QL + \frac{X}{P} L^{C+4} \right) \frac{X}{P} + (Q + XL^5)XL^{-C}\ .
\end{split}
\end{equation*}
Now, for every $S'(\sim PQ/d)$ and $S'(\sim 2PQ/d)$ we apply \eqref{5.21} (recall that $S'(\sim Q) \leqslant S(\sim Q,X)$) and sum over $d \leqslant L^C$ getting
\begin{equation*}
  \begin{split}
  X^{-1} S(\sim Q,X) \ll PQL^2 + (PQ)^{-1} X^{\frac{3}{2}} L^{11B + 2C + 47} + XL^{C-B} \\
  + P^{-2} X L^{C +4} + X L^{5-C}.
\end{split}
\end{equation*}
Here $B$, $C$, $P$ are arbitrary. We take $C = A + 6$, $B = C + A + 1 = 2A + 7$, so $11B + 2C + 47 = 24A + 136$ and $P= L^{28A + 147}$ getting
\begin{equation}
  S(\sim Q,X) \ll \left( Q L^{28A + 147} + Q^{-1} X^{\frac{3}{2}} L^{-4A - 11} + XL^{-A-1} \right) X .
  \label{6.3}
\end{equation}

If $Q>{\sqrt X}L^{-3A-10}$ then the middle term of \eqref{6.3} is covered by the last one.  
Finally, summing \eqref{6.3} over dyadic segments and incorporating \eqref{4.5}with $Q={\sqrt X} L^{-3A-10}$, we obtain the bound \eqref{3.9}. This completes the proof of the Theorem.

\end{document}